\documentclass[a4paper,12pt]{article}

\usepackage[dvipdfmx]{graphicx}
\usepackage{amsmath}
\usepackage{url}
\usepackage{txfonts}
\usepackage{multirow}
\usepackage[dvipdfmx,bookmarks=true,bookmarksnumbered=true,colorlinks=
true]{hyperref}

\begin{document}
\title{Practical Implementation of High-Order Multiple Precision Fully Implicit Runge-Kutta Methods with Step Size Control Using Embedded Formula}
\author{Tomonori Kouya\\Shizuoka Institute of Science and Technology\\2200-2 Toyosawa, Fukuroi, Shizuoka 437-8555 Japan\thanks{2010 Mathematics Subject Classification: 65L06, 65F10, 65G50}\thanks{Keywords and phrases: Implicit Runge-Kutta Method, Multiple Precision Floating-point Arithmetic, Iterative Refinement Method}}
\maketitle

\begin{quotation}\begin{center}\bf Abstract\end{center}
We propose a practical implementation of high-order fully implicit Runge-Kutta(IRK) methods in a multiple precision floating-point environment. Although implementations based on IRK methods in an IEEE754 double precision environment have been reported as RADAU5 developed by Hairer and SPARK3 developed by Jay, they support only 3-stage IRK families. More stages and higher-order IRK formulas must be adopted in order to decrease truncation errors, which become relatively larger than round-off errors in a multiple precision environment. We show that SPARK3 type reduction based on the so-called W-transformation is more effective than the RADAU5 type one for reduction in computational time of inner iteration of a high-order IRK process, and that the mixed precision iterative refinement method is very efficient in a multiple precision floating-point environment. Finally, we show that our implementation based on high-order IRK methods with embedded formulas can derive precise numerical solutions of some ordinary differential equations.
\end{quotation}

%
\section{Introduction}

Multiple precision floating-point (MP) arithmetic is an effective approach to solving ill-conditioned problems that cannot be solved precisely with IEEE754 double precison floating-point (DP) arithmetic. We have been developing BNCpack, a DP and MP numerical computation library based on MPFR\cite{mpfr}, an arbitrary precision and IEEE754 standard compatible floating-point arithmetic library, for the natural number arithmetic kernel in GNU MP (GMP)\cite{gmp} that is well-tuned for various CPU architectures. In this paper, we propose the implementation of a practical ordinary differential equation (ODE) solver based on BNCpack and high-order implicit Runge-Kutta (IRK) methods; its availability and efficiency are verified via numerical experiments.

DP ODE solvers based on IRK methods have been developed as RADAU5 (Radau IIA formula) and SPARK3 (selectable in Radau, Gauss, and Lobatto formulas) by Haier and Jay, respectively. Both IRK implementations support only 3-stage formulas, which is not sufficient to obtain precise numerical solutions in an MP environment. The use of MP floating-point arithmetic decreases round-off errors in an IRK process, thereby increasing truncation errors. More stages and higher-order IRK formulas are neccessary in an MP environment. Consequently, higher-dimensional nonlinear equations must be solved in high-order IRK processes. Such so-called inner iteration includes linear equations of the same dimension, which must be solved efficiently. In this process, RADAU5 reduces these coeffient matrices to complex diagonal matrices, and SPARK3 reduces them to real nonsymmetric tridiagonal matrices. In an MP environment, we must also accelerate these processes. In addition, we must be able to control the step sizes in an IRK process by using local error estimation at each discrete point in the given integration interval. For this reason, embedded formulas incidental to high-order IRK ones are neccessary. Althought Hairer proposed 4-stage embedded formula incidental to 3-stage Radau IIA formula, he did not describe explicitly the existence of embedded formulas for other IRK ones.

In this paper, we first state mathematical definitions and provide a framework for IRK algorithms.  In section 3, we explain the method of linear equations to be solved in the inner iteration of an IRK process. We show that the RADAU5 type reduction is not effective for high-order IRK methods, and that the mixed precision iterative refinement method can achive drastic acceleration, as shown via benchmark tests of a linear ODE. In section 4, we describe the derivation of embedded formulas for any IRK ones, and we discuss the A-stabilities. In section 5, we describe numerical experiments conducted for some test problems in order to demonstrate the high performance of our implementations. Finally, we conclude this paper and discuss the scope for future studies in section 6.

%
\section{Algorithm of Implicit Runge-Kutta Method}

We define the initial value proble (IVP) of the ODE to be solved as
\begin{equation}
\left\{\begin{array}{l}
	\displaystyle\frac{d\mathbf{y}}{dx} = \mathbf{f}(x, \mathbf{y}) \in\mathbb{R}^n \\
\\
	\mathbf{y}(x_0) = \mathbf{y}_0.
\end{array}\right. \label{eqn:ode}
\end{equation}
The integration interval is given as $[x_0, \alpha]$.

For this IVP of the ODE, we divide the integration interval into $x_0$, $x_1 := x_0 + h_0 $, ..., $x_{k+1} := x_k + h_k$, .... In order to obtain numerical solutions $\mathbf{y}_k \approx \mathbf{y}(x_k)$ in each step by using an $m$-stage IRK method, we must solve the following nonlinear system of equations (*). This process of solving with various iterative methods called inner iteration in IRK methods.
\begin{equation}
\begin{split}
& (*) \left\{\begin{array}{rcl}
	\mathbf{k}_1 &=& \mathbf{f}(x_k + c_1 h_k, \mathbf{y}_k + h_k\cdot \sum^{m}_{j=1} a_{1j}\mathbf{k}_j)\\
	\mathbf{k}_2 &=& \mathbf{f}(x_k + c_2 h_k, \mathbf{y}_k + h_k\cdot \sum^{m}_{j=1} a_{2j}\mathbf{k}_j)\\
	&\vdots& \\
	\mathbf{k}_m &=& \mathbf{f}(x_k + c_m h_k, \mathbf{y}_k + h_k\cdot \sum^{m}_{j=1} a_{mj}\mathbf{k}_j)\\
\end{array}\right. \\
\mathbf{y}_{k+1} &:= \mathbf{y}_k + h_k\cdot \sum^{m}_{j=1} b_j \mathbf{k}_j
\end{split} \label{eqn:eqn_irk}
\end{equation}
where the constant coefficients in the IRK formula, $c_1$, ..., $c_m$, $a_{11}$, ..., $a_{mm}$, $b_1$, ..., $b_m$, can be expressed as follows:
\[\begin{array}{c|cccc}
	c_1    & a_{11} & a_{12} & \cdots & a_{1m}   \\
	c_2    & a_{21} & a_{21} & \cdots & a_{2m}   \\
	\vdots & \vdots & \vdots &        & \vdots  \\
	c_m    & a_{m1} & a_{m2} & \cdots & a_{m,m}\\ \hline
	\      & b_1    & b_2    & \cdots & b_m \\
\end{array} = \begin{array}{c|c}
	\mathbf{c} & A \\ \hline
	           & \mathbf{b}^T
\end{array}.\]

All computations are the same in each discretization point $x_k$; hence, we consider only the computation $\mathbf{y}_0 \rightarrow \mathbf{y}_1 \approx \mathbf{y}(x_0 + h_0) = \mathbf{y}(x_0 + h)$.

\paragraph{Quasi-Newton Method}

If Newton method is adpoted as the numerical argorithm in inner iteration, the algorithm is as follows. The initial guesses are $\mathbf{k}^{(0)}_1$, ..., $\mathbf{k}^{(0)}_m$, and the approximations of unknowns $\mathbf{k}_1$, $\mathbf{k}_2$, ..., $\mathbf{k}_m$ are calculated by iterating the computations as
\[\left[\begin{array}{c}
	\mathbf{k}^{(l+1)}_1 \\
	\mathbf{k}^{(l+1)}_2 \\
	\vdots \\
	\mathbf{k}^{(l+1)}_m
\end{array}\right] := \left[\begin{array}{c}
	\mathbf{k}^{(l)}_1 \\
	\mathbf{k}^{(l)}_2 \\
	\vdots \\
	\mathbf{k}^{(l)}_m
\end{array}\right] - J^{-1}(\mathbf{k}^{(l)}_1, ..., \mathbf{k}^{(l)}_m) \left[\begin{array}{c}
	\mathbf{k}^{(l)}_1 - \mathbf{f}(x_0 + c_1 h, \mathbf{y}_0 + h \sum^{m}_{j=1} a_{1j} \mathbf{k}^{(l)}_j)\\
	\mathbf{k}^{(l)}_2 - \mathbf{f}(x_0 + c_2 h, \mathbf{y}_0 + h \sum^{m}_{j=1} a_{2j} \mathbf{k}^{(l)}_j)\\
	\vdots \\
	\mathbf{k}^{(l)}_m - \mathbf{f}(x_0 + c_m h, \mathbf{y}_0 + h \sum^{m}_{j=1} a_{mj} \mathbf{k}^{(l)}_j)
\end{array}\right] \]
where $J(\mathbf{k}^{(l)}_1, \mathbf{k}^{(l)}_2, ..., \mathbf{k}^{(l)}_m) \in\mathbb{R}^{mn\times mn}$ is given by
\[ J(\mathbf{k}^{(l)}_1, \mathbf{k}^{(l)}_2, ..., \mathbf{k}^{(l)}_m) = \left[\begin{array}{c|c|c|c}
	I_n - J_{11} & -J_{12}      & \cdots & -J_{1m} \\ \hline
	-J_{21}      & I_n - J_{22} & \cdots & -J_{2m} \\ \hline
	\vdots       & \vdots       &        & \vdots \\ \hline
	-J_{m1}      & -J_{m2}      & \cdots & I_n - J_{mm} \\
\end{array}\right] , \]
$I_n$ is an $n$-dimensional identity matrix, and $J_{pq}$ is given by
\[ J_{pq} = h a_{pq} \frac{\partial}{\partial \mathbf{y}} \mathbf{f}(x_0 + c_p h, \mathbf{y}_0 + h \sum^{m}_{j=1} a_{pj} \mathbf{k}^{(l)}_j)\in \mathbb{R}^{n\times n} .\]

In order to compute this part, we solve the following $mn$-dimensional system of linear equations with the coefficient matrix $J(\mathbf{k}^{(l)}_1, \mathbf{k}^{(l)}_2, ..., \mathbf{k}^{(l)}_m)$ for unknowns $[\overline{\mathbf{z}_1}\ \overline{\mathbf{z}_2}\ ...\ \overline{\mathbf{z}_n}]^T$.
\[ \begin{split}
& J(\mathbf{k}^{(l)}_1, \mathbf{k}^{(l)}_2, ..., \mathbf{k}^{(l)}_m) \left[\begin{array}{c}
	\overline{\mathbf{z}_1} \\
	\overline{\mathbf{z}_2} \\
	\vdots \\
	\overline{\mathbf{z}_m}
\end{array}\right] = \left[\begin{array}{c}
	\mathbf{k}^{(l)}_1 - \mathbf{f}(x_0 + c_1 h, \mathbf{y}_0 + h \sum^{m}_{j=1} a_{1j} \mathbf{k}^{(l)}_j) \\
	\mathbf{k}^{(l)}_2 - \mathbf{f}(x_0 + c_2 h, \mathbf{y}_0 + h \sum^{m}_{j=1} a_{2j} \mathbf{k}^{(l)}_j) \\
	\vdots \\
	\mathbf{k}^{(l)}_m - \mathbf{f}(x_0 + c_m h, \mathbf{y}_0 + h \sum^{m}_{j=1} a_{mj} \mathbf{k}^{(l)}_j)
\end{array} \right] \\
\end{split} \]
In this process, we need to have $O(m^2n^2)$ memory in order to store $J(\mathbf{k}^{(l)}_1, \mathbf{k}^{(l)}_2, ..., \mathbf{k}^{(l)}_m)$.

Actually, to accelerate this process, the following fixed coefficient matrix $J$ is used.

\[J(\mathbf{k}^{(l)}_1, \mathbf{k}^{(l)}_2, ..., \mathbf{k}^{(l)}_m) = J(\mathbf{k}^{(0)}_1, \mathbf{k}^{(0)}_2, ..., \mathbf{k}^{(0)}_m) \]

We call this method ``Quasi-Newton Method."

\paragraph{Simplified Newton Method}

For inner iteration in IRK methods, a more simplified Newton method is used conventionally\cite{hairer}. RADAU5 and SPARK3 select the simplified Newton method.

In this process, $J_{pq}$ used in the Quasi-Newton method is fixed as
\[J_{pq} := h a_{pq} \frac{\partial}{\partial \mathbf{y}} \mathbf{f}(x_0, \mathbf{y}_0) = h a_{pq} J ,\]
and then, we use $Y_i = \mathbf{y}_0 + h \sum^{m}_{j=1} a_{ij} \mathbf{f}(x_0 + c_i h, Y_j)$, which is an alternative to $\mathbf{k}^{(l)}_i$. Thus, we can express the system of linear equations to be solved as
\begin{equation}
	(I_m\otimes I_n - h A\otimes J) \mathbf{Z} = -\mathbf{F}(\mathbf{Y}) \in \mathbb{R}^{mn}, \label{eqn:simplified_newton_linear_eq}
\end{equation}
where
\[\mathbf{F}(\mathbf{Y}) = \left[\begin{array}{c}
	Y_1 - \mathbf{y}_0 - h\sum^m_{j=1} a_{1j}\mathbf{f}(x_0 + c_1 h, Y_1) \\
	\vdots\\
	Y_m - \mathbf{y}_0 - h\sum^m_{j=1} a_{mj}\mathbf{f}(x_0 + c_m h, Y_m)
\end{array}\right].\]

\section{Acceleration of Inner Iteration in IRK process}

In this section, we treat the simplified Newton method as inner iteration. In this case, the system of linear equations to be solved at each step of inner iteration is expressed in (\ref{eqn:simplified_newton_linear_eq}). We can decrease the number of computations by applying reductions based on similarity transformations to IRK matrix $A$. RADAU5 employs complex diagonalization, and SPARK3 employs real unsymmetric tridiagonalization, which is called W-transformation. The latter transformation is better than the former because it does not treat complex arithmetic, and it can avoid the ill-conditioned transformations. Therefore, we employ SPARK3 type reduction for the implementation of high-order IRK methods.  In the rest of this section, we compare RADAU5 type and SPARK3 type reductions, and we show that the acceleration due to the mixed precision iterative refinement method can drastically reduce the computational time of an IRK process.

\subsection{Comparison between RADAU5 type and SPARK3 type Reductions}

RADAU5 type reduction\cite{hairer} is based on the fact that IRK matrix $A$ can be transformed into a complex diagonal matrix as
\[ \Lambda = \mbox{diag}(\lambda_1, \cdots, \lambda_m) = S A S^{-1} .\]
In these cases of IRK Radau, Gauss, and Lobatto formula families, eigenvalues $\lambda_i$ are generally complex numbers.

If RADAU5 type reduction is employed in the simplified Newton method, the system of linear equations to be solved has the coefficient matrix obtained by similarity transformation with $S\otimes I_n$ and $S^{-1}\otimes I_n$,  as
\[ \begin{split}
	&(S\otimes I_n)(I_m\otimes I_n - hA \otimes J)(S^{-1}\otimes I_n) = I_m\otimes I_n - h\Lambda \otimes J \\
	& = \left[\begin{array}{ccc}
		I_n - h\lambda_1 J & & \\
		& \ddots & \\
		& & I_n - h\lambda_m J
\end{array}\right].
\end{split}\]

\begin{figure}[htpb]
\begin{center}
\includegraphics[width=.5\textwidth]{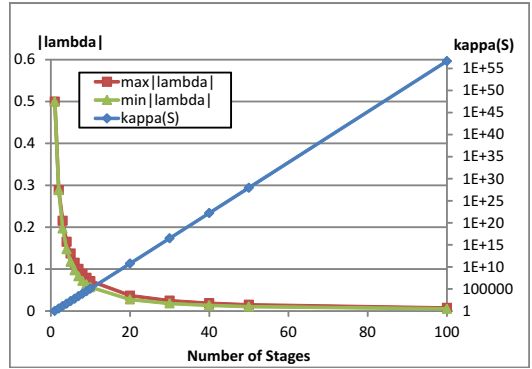}
\caption{$\kappa_2(S)$ and absolute values of eigenvalues of $A$ for IRK Gauss formulas}\label{fig:complex_diagonal}
\end{center}
\end{figure}

RADAU5 type reduction has two advantages: 1. The inner iteration is completely parallelizable, and 2. the order of the required memory is $O(2mn)$. However, $\Lambda, S$, and $S^{-1}$ are complex matrices, and $\kappa_2(S) = \|S\|_2 \|S^{-1}\|_2\rightarrow$ $\infty$ $(m\rightarrow \infty)$; hence, round-off errors increase in inner iteration, especially for higher order IRK formulas,  as shown in \figurename\ \ref{fig:complex_diagonal}.

On the other hand, SPARK3 type reduction is called block triangulation through W-transformation\cite{hairer}, given by 
\[\begin{split}
	X &= W^T BA W = \left[\begin{array}{ccccc}
		1/2     & -\zeta_1 &             &              & \\
		\zeta_1 & 0        & \ddots      &              & \\
		        & \ddots   & \ddots      & -\zeta_{m-2} & \\
		        &          & \zeta_{m-2} & 0            & -\zeta_{m-1} \\
	            &          &             & \zeta_{m-1}  & 0
\end{array}\right], \\
	\mbox{where }& \\
 	w_{ij} = \tilde{P}_{j-1}(c_i) &= \sqrt{2(j-1)+1} \sum^{j-1}_{k=1} (-1)^{j+k-1} \left(\begin{array}{c}
	j-1 \\
	k
\end{array}\right) \left(\begin{array}{c}
	j+k-1 \\
	k
\end{array}\right) c_i^k \\
	\mbox{$\tilde{P}_{j-1}(x)$}&: \mbox{$j-1$-th shifted Legendre polynomial} \\
	\zeta_i &= \left(2\sqrt{4i^2-1}\right)^{-1},\ B = \mbox{diag}(\mathbf{b}),\ D = W^T B W = \mbox{diag}(1\ 1\ \cdots\ 1).
\end{split} \]

By using W-transformation, we can obtain the reduced coefficient matrix in the simplified Newton method \cite{jay99} as
\[\begin{split}
& (W^TB\otimes I_n)(I_m\otimes I_n - h A\otimes J)(W\otimes I_n) \\
&= D\otimes I_n - hX\otimes J = \left[\begin{array}{ccccc}
	E_1 & F_1    &         & & \\
	G_1 & E_2    & F_2     & & \\
	    & \ddots & \ddots  & \ddots  & \\
		&        & G_{m-2} & E_{m-1} & F_{m-1} \\
	    &        &         & G_{m-1} & E_m
\end{array}\right]
\end{split}, \]
where 
\[ \begin{split}
	E_1 &= I_n - \frac{1}{2}hJ,\ E_2 = \cdots = E_s = I_n \\
	F_i & = h\zeta_iJ,\ G_i = -h\zeta_iJ\ (i = 1, 2, ..., m-1)
\end{split}.\]

In the implementation of SPARK3, the left preconditioner matrix $P$, 
\[ P = \left[\begin{array}{ccccc}
	\tilde{E}_1 & F_1    &         & & \\
	G_1 & \tilde{E}_2    & F_2     & & \\
	    & \ddots & \ddots  & \ddots  & \\
		&        & G_{m-2} & \tilde{E}_{m-1} & F_{m-1} \\
	    &        &         & G_{m-1} & \tilde{E}_m
\end{array}\right] \approx D\otimes I_n - hX\otimes J \]
is prepared, and then,  the precondtioned system of linear equation 
\[ P^{-1} (D\otimes I_n - hX\otimes J) (W\otimes I_n)^{-1} \mathbf{Z}= P^{-1}(W^TB\otimes I_n)(-\mathbf{F}(\mathbf{Y})) \]
is set to be solved for $\mathbf{Z}$. Jay maintained that the left precondition can decrease the number of iterations and accelerate Richardson iteration and GMRES (Generalized Minimal RESidual) methods\cite{jay99}. However, our numerical experiments show that  such preconditioning  increases the computational time in MP environment; hence, our current implementation employs such preconditions in DP environment, and does not employ in MP environment.

In comparison with RADAU5 type reduction, SPARK3 type reduction is better because all computations constitute real number arithmetic, and the order of the memory required for the coefficient matrix is $O(3mn)$. In addition, the similarity transformation matrix $W$ can remain well-conditioned when the number of stages $m$ is large, and hence, the effect of round-off errors occuring in the similarity transformation is small. \tablename\ \ref{table:cond_w} shows that the condition number of $\kappa_\infty(W)=\|W\|_\infty\|W^{-1}\|_\infty$ is much smaller than $\kappa_\infty(S)$ used in RADAU5 type reduction.

\begin{table}
\begin{center}
\caption{Condition numbers of two kinds of similarity transformation matrices}\label{table:cond_w}
\begin{tabular}{|c|c|c|c|c|c|c|}\hline
$m$                & 3      & 5      & 10                & 15               & 20                   & 50 \\ \hline\hline
$\kappa_\infty(S)$ & $22.0$ & $388$  & $3.28\times 10^5$ & $2.81\times 10^{8} $ & $2.11\times 10^{11}$ & $4.25\times 10^{28}$ \\ \hline
$\kappa_\infty(W)$ & $3.24$ & $6.27$ & $16.4$            & $29.3$ & $44.5$               & $172$ \\ \hline
\end{tabular}
\end{center}
\end{table}

%
\subsection{Acceleration by using Mixed Precision Iterative Refinement Method}

The mixed precision iterative refinement method was originally proposed by Moler in 1967\cite{moler_iterative_ref}, and then Buttari et al. showed that their revised algorithm exhibits high performance in many current computing environments. If $S$-digit floating-point arithmetic can be executed more efficiently than $L$ $(>>S)$-digit arithmetic, the system of linear equations
\[ \mathbf{f}(\mathbf{x}) = C\mathbf{x} - \mathbf{d} \]
can be solved by Newton method and an appropriate linear solver, shown by the following algorithm:
\begin{align}
\mbox{Solve}\ & \ C^{[S]}\mathbf{x}_0^{[S]} = \mathbf{d}^{[S]}\ \mbox{for}\ \mathbf{x}_0^{[S]}. \label{eqn:simple_iterative_ref_l1} \\
\mathbf{x}_0^{[L]} &:= \mathbf{x}_0^{[S]} \nonumber \\
\mbox{For}\ & \ k=0, 1, 2, ... \nonumber \\
	& \mathbf{r}_k^{[L]} := \mathbf{d}^{[L]} - C^{[L]}\mathbf{x}_k^{[L]} \label{eqn:simple_iterative_ref1} \\
	& \mathbf{r}_k'^{[L]} := \mathbf{r}_k^{[L]} / \|\mathbf{r}_k^{[L]}\| \nonumber \\
	& \mathbf{r}_k'^{[S]} := \mathbf{r}_k'^{[L]} \nonumber \\
	& \mbox{Solve}\ C^{[S]}\mathbf{z}_{k}^{[S]} = \mathbf{r}_k'^{[S]}\ \mbox{for}\ \mathbf{z}_k^{[S]}. \label{eqn:simple_iterative_ref_l2} \\
	& \mathbf{z}_k^{[L]} := \mathbf{z}_k^{[S]} \nonumber \\
	& \mathbf{x}_{k+1}^{[L]} := \mathbf{x}_k^{[L]} + \|\mathbf{r}_k^{[L]}\| \mathbf{z}_k^{[L]} \label{eqn:simple_iterative_ref3} \\
	& \mbox{Check convergence of $\mathbf{x}_{k+1}$.} \nonumber
\end{align}
where $[S]$ and $[L]$ denote the values expressed and computed in $S$- and $L$-digit floating-point arithmetic, respectively. The above algorithm is the $S$-$L$ mixed precision iterative refinement method for a system of linear equations. The part of  (\ref{eqn:simple_iterative_ref_l1}) and (\ref{eqn:simple_iterative_ref_l2}) theoretically give the solution, and hence, we do not need the above iteration. However, we cannot obtain the true solutions owing to the use of finite precision floating-point arithmetic in these parts. Thus, the residuals $\mathbf{r}_k$ are not all zero. The above algorithm executes some iterations so that the residuals $\mathbf{r}_k$ tend to zero. The number of digits required to compute the residuals increases with the precision of the obtained approximation $\mathbf{x}_k$.

Buttari et. al. also proved that the sufficient condition for convergence is satisfied if the precision of the computation (\ref{eqn:simple_iterative_ref_l2}) is less than that of the residual (\ref{eqn:simple_iterative_ref1}) when the condition number $\kappa(A)=\|A\|\|A^{-1}\|$ is smaller than the precision in the computing environment. In addition, they showed that the combination of (\ref{eqn:simple_iterative_ref1}) and (\ref{eqn:simple_iterative_ref3}) computed DP arithmetic, and (\ref{eqn:simple_iterative_ref_l1}) and (\ref{eqn:simple_iterative_ref_l2}) computed in SP(IEEE754 single-precision) arithmetic can accelerate the solution of relatively less ill-conditioned systems of linear equations via benthmark\cite{mixed_prec_iterative_ref}\cite{lawn175}.

Their mixed precision iterative refinement method is just fitted to solve the system of linear equations (\ref{eqn:simplified_newton_linear_eq}) in inner iteration of IRK methods. This tends to be well-conditioned if the step size $h_k$ becomes smaller, and hence, the sufficient condition of convergence is satisfied in almost cases in an IRK process. In particular, for MP arithmetic, the application to the DP-MP type mixed precision iterative refinement method in which (\ref{eqn:simple_iterative_ref_l2}) is computed by DP arithmetic can drastically accelerate the entire computational time of an IRK process. We support direct methods and Krylov subspace methods in our current implementation of mixed precision iterative refinement method.

In the rest of this section, we show the acceleration of an IRK method with DP-MP type mixed precision iterative refinement method via benchmark tests applied to a linear ODE.

%

Our test problem is constructed by using a real normal matrix $R$ comprising uniform random numbers, and by its inverse matrix $R^{-1}$, as follows:

\[\begin{split}
& \left\{\begin{array}{l}
	\displaystyle\frac{d\mathbf{y}}{dx} = -(R\ \mbox{diag}(n, n-1, ..., 1)\ R^{-1})\ \mathbf{y}\ \in\mathbb{R}^{128} \\
	\\
	\mathbf{y}(0) = [1\ ...\ 1]^T \\
\end{array}\right. \\
& \mbox{Integration Interval:} [0, 20].
\end{split}\]
The precision of MP arithmetic is fixed at $50$ decimal digits (167 bits). We compute $\mathbf{y}_1$ by using the $m$-stage $2m$ order Gauss formulas ($m=3,4, ..., 12$ and $h=1/2$), and we compare the following 4 algorithms:
\begin{enumerate}
\item Quasi-Newton method with the DP-MP type mixed precision iterative refinement method based on the direct method and without reduction: ``Iter.Ref-DM"
\item Simplified Newton method with the simple direct method and with SPARK3 type reduction: ``W-Trans"
\item Simplified Newton method with the MP (25 decimal digits)-MP iterative refinement method with SPARK3 type reduction: ``W-Iter.Ref-MM"
\item Simplified Newton method with the DP-MP method with SPARK3 type reduction: ``W-Iter.Ref-DM"
\end{enumerate}

\begin{figure}[htbp]
\begin{center}
\includegraphics[width=.8\textwidth]{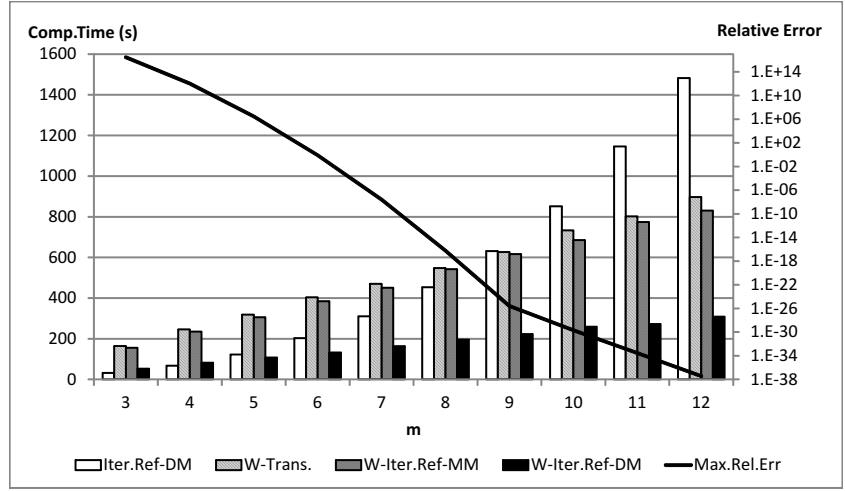}
\caption{Performance of IRK methods: in case of 128-dimensional and 50 decimal digits computation}\label{fig:irk_jay_bench}
\end{center}
\end{figure}

All computations are executed on an Intel Core i7 920 + CentOS 5.4 x86\_64 machine with gcc 4.1.2 + BNCpack 0.8 + MPFR 3.1.0/GMP 5.0.2. For convenient comparison, the maximum relative error in approximation $\mathbf{y}_1$ is expressed as a line graph in \figurename\ \ref{fig:irk_jay_bench}. All 4 algorithms can obtain the same accuracy of $\mathbf{y}_1$.

Although the relative errors over 9 stages undergo less reduction because of the effect of round-off errors, all numerical results from 3- to 12-stages IRK Gauss formulas are precise.

The computational time increases with number of stages employed in the IRK formulas. The Quasi-Newton method without reduction in inner iteration is competitive with less than four stages, but is extremely slow with over 9 stages. On the other hand, the computational times of the three algorithms with SPARK3 type reduction is proportional to the number of stages. Moreover, the application of the DP-MP type iterative refinement method can drastically accelerate IRK processes, especially for formulas with over 5 stages formulas. The maximum obtained speedup ratio is 4.8.

%
\section{Derivation of Embedded Formula and Step Size Selection}

In order to implement actual ODE solvers, we need a mechanism of step size selection, based on the local error estimation at each discretized point. In the case of Runge-Kutta (RK) methods, embedded formulas that can reduce the amount of computation and the number of integrated function calls are widely used. For explicit RK (ERK) methods, many ODE solvers have been developed on the basis of the embedded formulas proposed by Fehlberg and Dormand-Prince. On the other hand, the embedded formula incident to 3-stage Radau IIA formula, which is proposed by Hairer\cite{hairer2} is only one for IRK methods, and it is constructed by a combination of the original formula and a lower order one. Hairer suggested the existence of the same type embedded formulas for other IRK formulas. In this section, we show how to automatically derive embedded formulas for IRK Gauss formulas which have wider A-stability regions.

%
\subsection{Hairer's Embedded Formula}

Hairer's embedded formula for the 3-stage 5th order Radau IIA formula is used in RADAU5,  and it can be expressed as follows. The original Radau IIA is 
\[\begin{array}{c|ccc}
\frac{4-\sqrt{6}}{10} & \frac{88-7\sqrt{6}}{360} & \frac{296-169\sqrt{6}}{1800} & \frac{-2+3\sqrt{6}}{225} \\
\frac{4+\sqrt{6}}{10} & \frac{296+169\sqrt{6}}{1800} & \frac{88+7\sqrt{6}}{360} & \frac{-2-3\sqrt{6}}{225} \\
1                     & \frac{16-\sqrt{6}}{36} & \frac{16+\sqrt{6}}{36} & \frac{1}{9} \\ \hline
           & \frac{16-\sqrt{6}}{36} & \frac{16+\sqrt{6}}{36} & \frac{1}{9} \\
\end{array} \]
and its embedded formula is the 4-stage formula
\begin{equation}
\begin{array}{c|cc}
c_0 = 0    & 0          & \mathbf{0}^T \\
\mathbf{c} & \mathbf{0} & A \\ \hline
           & \gamma_0   & \hat{\mathbf{b}}^T
\end{array} \label{eqn:irk_embed}
\end{equation}
where $\gamma_0$ is any non-zero constant, as recommended by Hairer for the real eigenvalue of IRK matrix $A$ in order to reduce the number of computations in local error estimation. Moreover, $\hat{\mathbf{b}}$ is expressed as 
\begin{equation}
	\hat{\mathbf{b}} = \left[\begin{array}{c}
		\hat{b}_1 \\
		\hat{b}_2 \\
		\hat{b}_3
\end{array}\right] = \left[\begin{array}{c}
		b_1 - \frac{2+3\sqrt{6}}{6}\gamma_0 \\
		b_2 - \frac{2-3\sqrt{6}}{6}\gamma_0 \\
		b_3 - \frac{\gamma_0}{3}
	\end{array}\right].
\end{equation}
This $\hat{\mathbf{b}}$ satisfies the simplifying assumption $B(3)$\ \cite{hairer}
\[ \mbox{$B(3)$:\ } c_0^{q-1} \gamma_0 + \sum^3_{i=1} \hat{b}_i c_{i}^{q-1} = 1/q\ (q=1, 2, 3). \]
The other coefficients $\mathbf{c}$ and $A$ are the same as the original formula, and hence, the simplifying $C(3)$, given by 
\[ \mbox{$C(3)$:\ } \sum^3_{j=1} a_{ij} c_j^{q-1} = c_i^q / q\ (i=1, 2, 3, q=1, 2, 3), \]
is automatically satisfied. Therefore, the given embedded formula (\ref{eqn:irk_embed}) is of 3rd order at least.

%
\subsection{Derivation of Embedded Formulas for any IRK ones}

Hairer's embedded formulation can be generally extended to other IRK formulas. In order to satisfy the simplifying assumption $B(m)$, $\hat{\mathbf{b}}$ is the solution of the following system of linear equations with a Vandermonde matrix of coefficients,
{\small \begin{equation}
\left[\begin{array}{cccc}
	1 & 1 & \cdots & 1 \\
	c_1 & c_2 & \cdots & c_m \\
	\vdots & \vdots & & \vdots \\
	c_1^{m-1} & c_2^{m-1} & \cdots & c_m^{m-1}
\end{array}\right] \left[\begin{array}{c}
	\hat{b}_1 \\
	\hat{b}_2 \\
	\vdots \\
	\hat{b}_m
\end{array}\right] = \left[\begin{array}{c}
	1 - \gamma_0 \\
	1/2 \\
	\vdots \\
	1/m
\end{array}\right] .
\label{eqn:irk_embed_intro}
\end{equation}}
By solving the above equation, we can obtain the approximation $\hat{\mathbf{y}}_{k+1}$ derived by the $m$-th order embedded formula,
\begin{equation}
	\hat{\mathbf{y}}_{k+1} := \mathbf{y}_k + h_k\left(\gamma_0 \mathbf{f}(x_k, \mathbf{y}_k) + \sum^{m}_{j=1} \hat{b}_j Y_j\right). \label{eqn:embed_formula}
\end{equation}

\begin{figure}[htpb]
\begin{center}
\includegraphics[width=.65\textwidth]{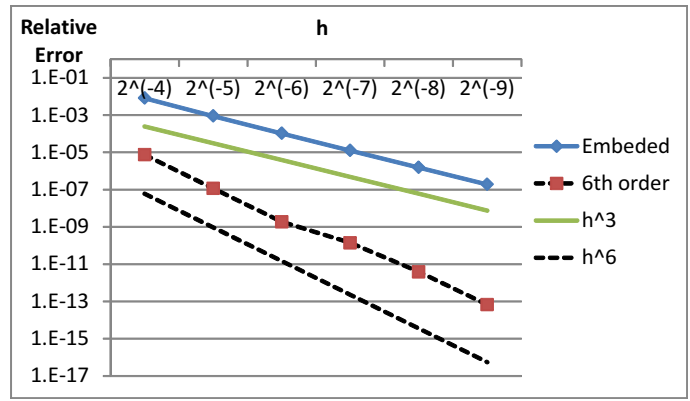}
\caption{Relative error of the 3rd order embedded formula derived from 3-stage 6th order IRK formula}\label{fig:mxy}
\end{center}
\end{figure}

When we fix $\gamma_0 = 1/8$ for the reason described in the next section, we compute the relative errors of the approximation $\hat{y}_k$ in (\ref{eqn:mxy}) at $x = 10$ by using the embedded formula derived from the 3-stage 6th order Gauss IRK formula.
\begin{equation}
\begin{array}{l}
	\displaystyle\frac{dy}{dx} = -xy \\
	y(0) = 1 \\
\mbox{Integration Interval:} [0, 10]
\end{array} \label{eqn:mxy}
\end{equation}
As a result, we can confirm that the order of the embedded formula is 3, as shown in \figurename\ \ref{fig:mxy}. For comparison, the original 6-th order $y_k$ is also plotted in the same figure.

In our implementation of the MP ODE solver based on high-order IRK Gauss formulas, we currently use the following value $\|\mathbf{err}_k\|$ for local error estimation\cite{hairer2}.

\[ \|\mbox{\bf err}_k\| = \sqrt{\frac{1}{n}\sum^{n}_{j=1}\left(\frac{\left|\hat{y}^{(k+1)}_j - y^{(k+1)}_j\right|}{ATOL + RTOL\cdot \max\left(\left|y^{(k+1)}_j\right|, \left|y^{(k)}_j\right|\right)}\right)^2} \]

The next step size $h_{k+1}$ at $x_{k+1}$ is set as
\[ h_{k+1} := 0.9 \|\mathbf{err}_k\|^{1/(m+1)} h_k. \]

%
\subsection{A-stability Regions of Embedded Formulas}

\begin{figure}[htpb]
\begin{center}
\includegraphics[width=.4\textwidth]{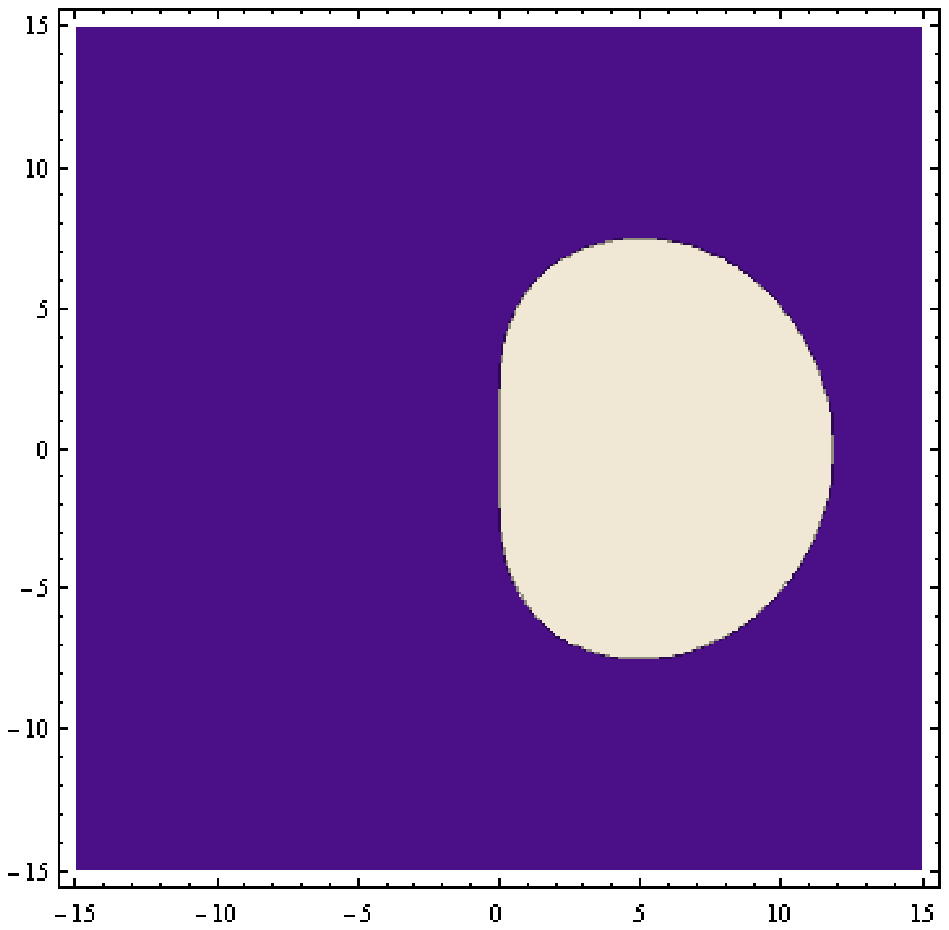}
\includegraphics[width=.4\textwidth]{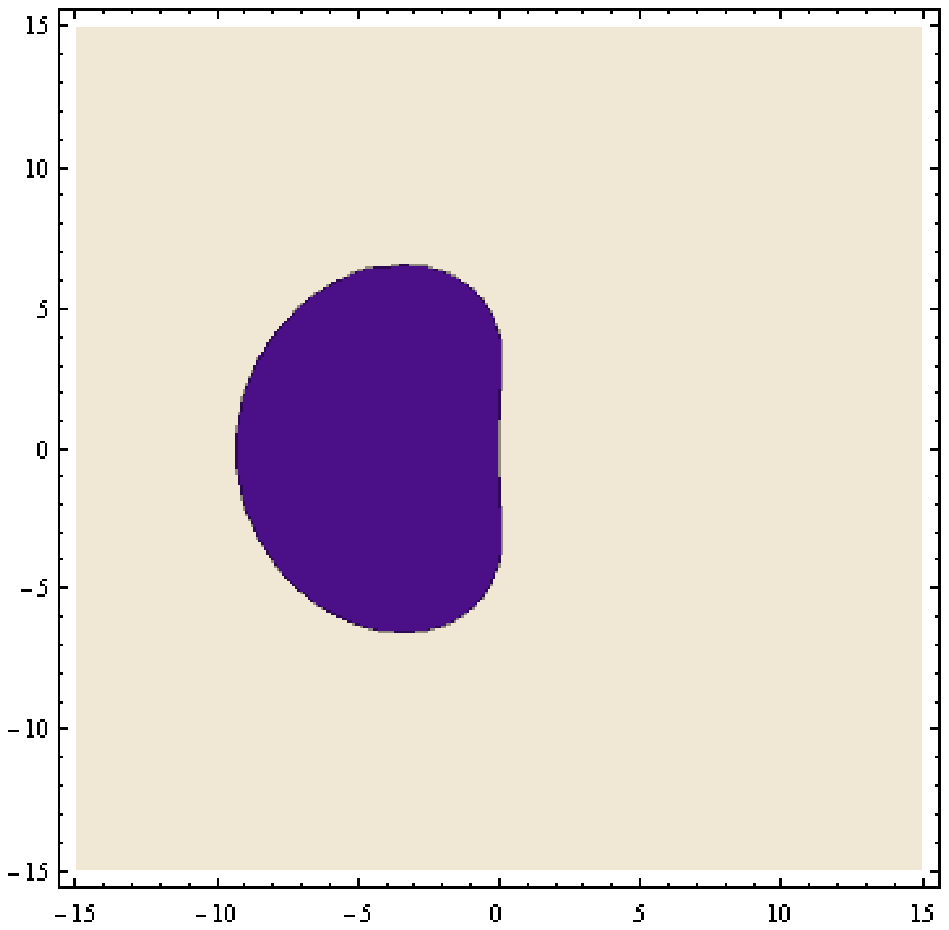}
\caption{Stability regions: 3-stage 5th order Radau IIA formula (left) and the corresponding Hairer's embedded formula (right).}\label{fig:irk_radau2}
\end{center}
\end{figure}

One advantage of IRK formulas can be A-stable; however, the embedded formulas derived from the original IRK formulas are not A-stable. Actually, the stability region of Hairer's embedded formula is narrower than that of the original Radau IIA formula, as shown in \figurename\ \ref{fig:irk_radau2}.

\begin{figure}[htpb]
\begin{center}
\includegraphics[width=.4\textwidth]{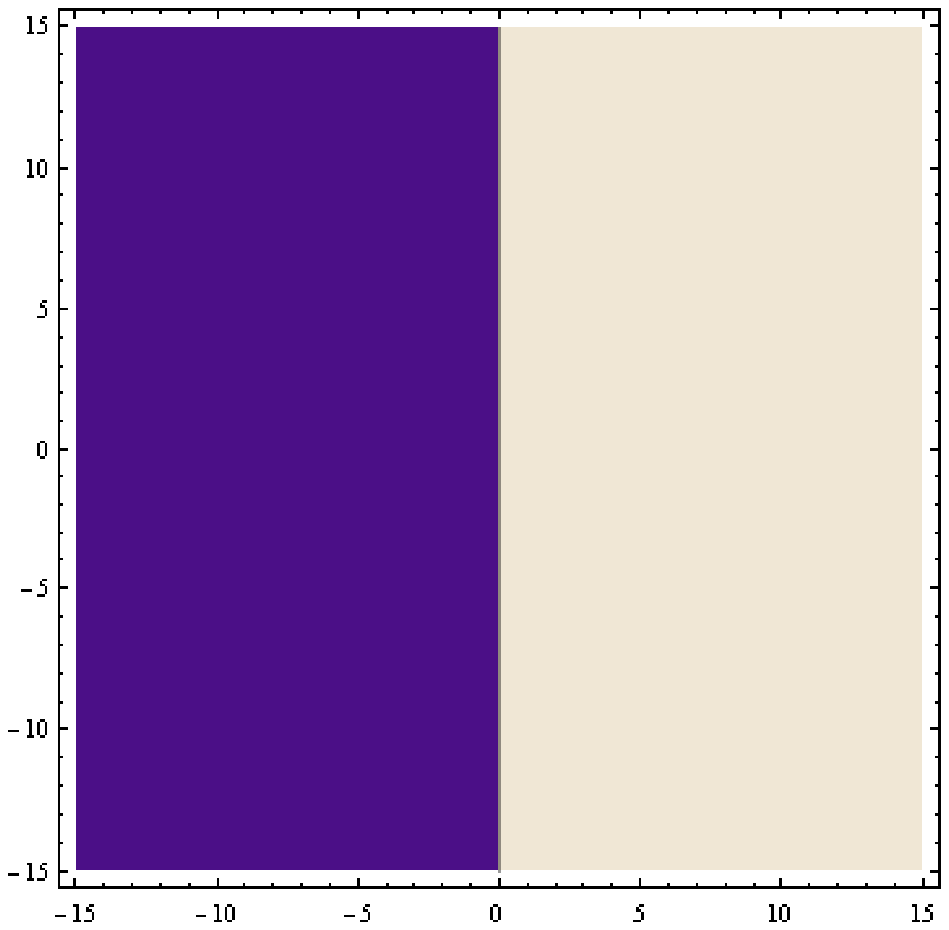}
\includegraphics[width=.4\textwidth]{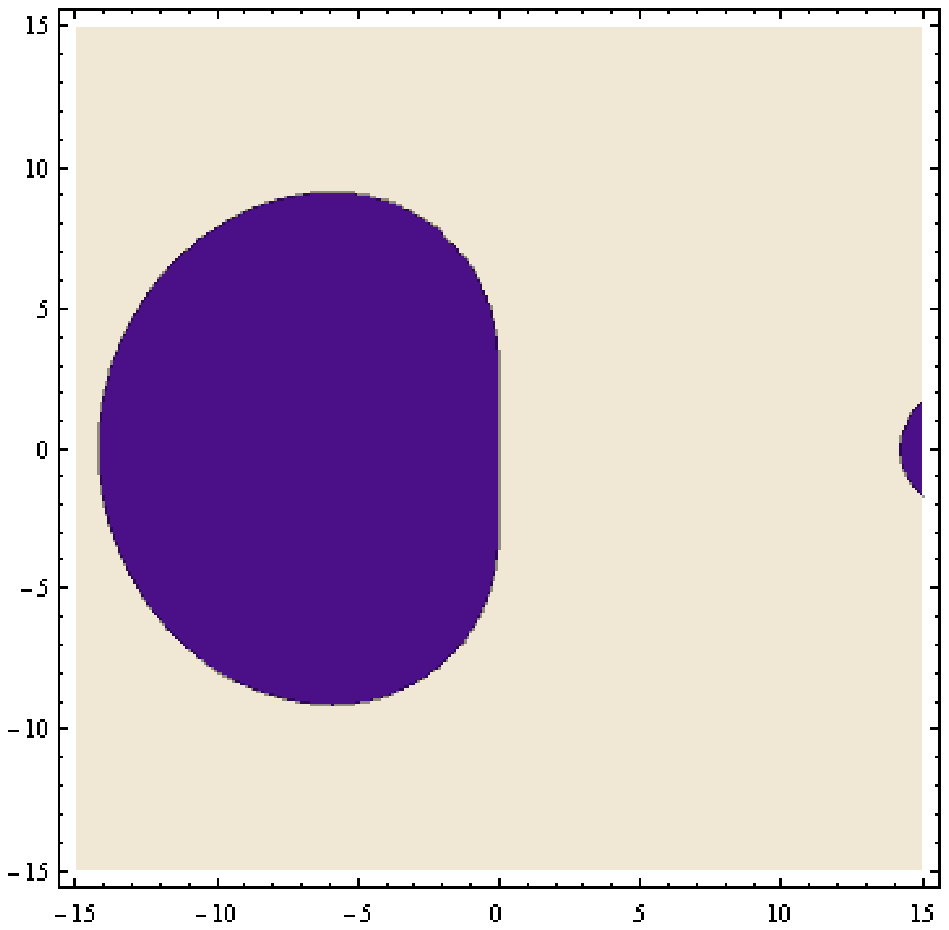}
\caption{Stability regions: 3-stage 6th order Gauss formula (left) and our embedded formula(left, $\gamma_0 = 1/8$)}\label{fig:irk_gauss3}
\end{center}
\end{figure}

We must select the parameter $\gamma_0 \not= 0$ at which the embedded formula can have a wider region. As a result, we currently consider that $\gamma_0 = 1/8$ is better because it is not too small and it can be expressed in powers of 2. In addition, its region is wider than that of Hairer's embedded formula. The stability region is shown in \figurename\ \ref{fig:irk_gauss3}.

%
\section{Numerical Experiments}

As described in previous sections, the following three techiques are applied to our MP ODE solver based on high-order IRK formulas:
\begin{enumerate}
\item Simplified Newton Method with SPARK3 type reduction in inner iteration
\item DP-MP type mixed precision iterative refinement method to accelarate inner iteration
\item Step size selection based on embedded formulas automatically derived
\end{enumerate}

In addition, users can select DP and MP Krylov subspace methods supporting a banded Jacobi matrix in the same way as SPARK3.

To evaluate the peformance of our implementation, we present the results of numerical experiments. All computations are executed on a Intel Core i7 920 + CentOS 5.4 x86\_64 machine (gcc 4.1.2 + BNCpack 0.8 + MPFR 3.1.0/GMP 5.0.2).

%
\subsection{Non Stiff Problem}

The Lorenz problem (\ref{eqn:lorenz}) is a well-known problem in complex systems, and it is not stiff; however, the accuracy of approximation is worse in a longer integration interval. Thus, we must use MP arithmetic in propotion to the length of the integration interval.

\begin{equation}
\begin{array}{l}
\left\{\begin{array}{ccl}
	\displaystyle\frac{d\mathbf{y}}{dx} &=& \left[\begin{array}{c}
		\sigma (- y_1 + y_2) \\
		-y_1 y_3 + ry_1 - y_2 \\
		y_1 y_2 - by_3
	\end{array}\right] \\
	\mathbf{y}(0) &=& [ 0\ 1\ 0 ]^T
\end{array}\right. \\
\mbox{Integration Interval:} [0, 50]
\end{array}, \label{eqn:lorenz}
\end{equation}
where $\sigma = 10$, $r = 470/19$, and $b = 8/3$. In the case of the above integration interval,  we lose around 13 decimal digits of the approximation of $\mathbf{y}(50)$. Hence, we select 70 decimal digits (233 bits), 10-stage 20th order and 15-stage 30th order Gauss formulas. The numerical results are shown in \tablename\ \ref{table:lorenz}.

\begin{table}
\begin{center}
\caption{Lorenz Problem}\label{table:lorenz}
\begin{tabular}{|c|c|c|c|c|}\hline
                 & \multicolumn{2}{|c|}{$RTOL=10^{-30}$, $ATOL=0$} & \multicolumn{2}{|c|}{$RTOL=10^{-50}$, $ATOL=0$} \\ \hline
                 & 10 stages & 15 stages & 10 stages & 15 stages \\ \hline
\# steps     & 41137    & 5112        & 2709021 & 91169 \\ \hline
Comp.Time (s)      & 192.0    & 64.3        & 12911.1 & 1112.2 \\ \hline
Max.Rel.Error & $3.9\times 10^{-19}$   & $4.4\times 10^{-19}$ & $3.8\times 10^{-39}$ & $4.9\times 10^{-39}$ \\
Min.Rel.Error & $7.3\times 10^{-21}$   & $8.3\times 10^{-21}$ & $7.1\times 10^{-41}$ & $9.1\times 10^{-41}$ \\ \hline
\end{tabular}
\end{center}
\end{table}

The numerical results using both 10-stage and 15-stage formulas indicate that we can obtain the appropriate accuracy of corresponding $RTOL$s.

%
\subsection{Stiff Problem}

\begin{figure}[htpb]
\begin{center}
\includegraphics[width=.9\textwidth]{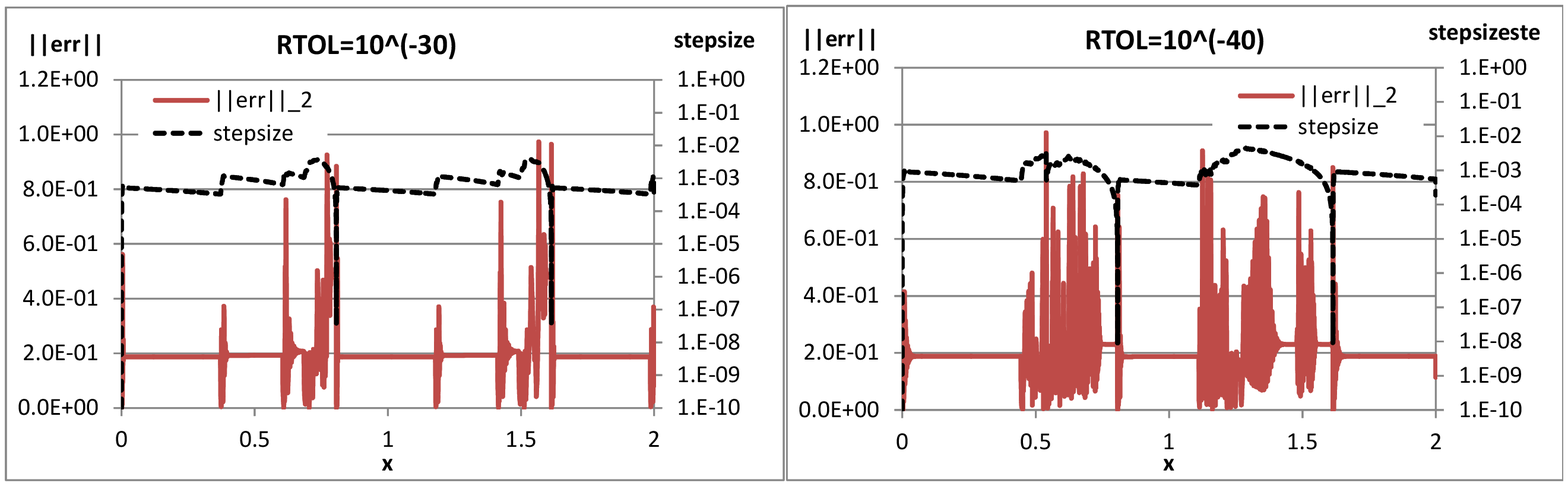}
\caption{History of local error estimation and step size: van der Pol equaiton}\label{fig:local_error_est}
\end{center}
\end{figure}

Next, we solve the van del Pol equation (\ref{eqn:vdpol}),  which is provided in Testset\cite{testset_ode}.
\begin{equation}
\begin{array}{l}
\left\{\begin{array}{ccl}
	\displaystyle\frac{d\mathbf{y}}{dx} &=& \left[\begin{array}{c}
		y_2 \\
		( (1 - y_1^2)y_2 - y_1 ) / 10^{-6}
	\end{array}\right] \\
	\mathbf{y}(0) &=& \left[\begin{array}{c}
		2 \\
		0
	\end{array}\right]
\end{array}\right. \\
\mbox{Integration Interval:} [0, 2]
\end{array} \label{eqn:vdpol}
\end{equation}
By using a 15-stage 30th order Gauss formula computed in 50 decimal digits (167 bits) MP arithmetic, we set $ATOL=0$ and $RTOL=10^{-40}$ or $10^{-30}$. In this case, \figurename\ \ref{fig:local_error_est} shows the history of $\|\mathbf{err}_k\|$ and step size.

As a result, we can obtain the appropriate approximations corresponding to each $RTOL$s listed in \tablename\ \ref{table:van_del_pol}.

\begin{table}
\begin{center}
\caption{Van del Pol equation: The number of steps, computational time and relative errors.}\label{table:van_del_pol}
\begin{tabular}{|c|c|c|}\hline
         & $RTOL=10^{-30}$ & $RTOL=10^{-40}$ \\ \hline
\# steps & 4325  & 6202 \\ \hline
Comp.Time (s)  & 127.8 & 208.9 \\ \hline
Max.Rel.Error & $1.2\times 10^{-29}$ & $1.0\times 10^{-39}$ \\
Min.Rel.Error & $2.2\times 10^{-36}$ & $8.6\times 10^{-47}$ \\ \hline
\end{tabular}
\end{center}
\end{table}

%
\subsection{1-Dimensional Brusselator Problem}

As a large-scale problem, we solve the 1-dimensional Brusselator problem\cite{spark3} by using our MP ODE solver.
\begin{equation}
\left\{\begin{array}{l}
	\frac{\partial u}{\partial t} = 1 + u^2 v - 4 + 0.02\cdot \frac{\partial^2 u}{\partial x^2} \\
	\frac{\partial v}{\partial t} = 3u - u^2v + 0.02\cdot \frac{\partial^2 v}{\partial x^2}
\end{array}\right. \label{eqn:bruss1d_org}
\end{equation}

The above original partial differential equation (\ref{eqn:bruss1d_org}) can be discretetized as a large-scale ODE as follows:
\begin{equation}
\begin{array}{l}
\left\{\begin{array}{l}
	\frac{d u_i}{d t} = 1 + u_i^2 v_i - 4 + 0.02\cdot \frac{u_{i+1} - 2u_i + u_{i-1}}{(\Delta x)^2} \\
	\frac{d v_i}{d t} = 3u_i - u_i^2v_i + 0.02\cdot \frac{v_{i+1} - 2v_i + v_{i-1}}{(\Delta x)^2} \\
u_0(t) = u_{N+1}(t) = 1,\ v_0(t) = v_{N+1}(t) = 3,\\
u_i(0) = 1 + \sin(2\pi i\Delta x),\ v_i(0) = 3
\end{array}\right. \ (i =1, 2, ..., N) \\
\mbox{Integration Interval:} [0, 10]
\end{array}
\label{eqn:bruss1d}
\end{equation}

We solve the above ODE for the following situation:
\begin{description}
\item[\bf Parameters] $N=500, n=2N=1000$, $\Delta x = 1/(N+1) = 1/501$
\item[\bf TOLs] $RTOL = ATOL = 10^{-30}$
\item[\bf Used linear solver] DP-MP($L=50$) type mixed precision BiCGSTAB method with band matrix-vector multiplication.
\end{description}

In many large-scale problems, $J$ can be sometimes expressed as a sparse matrix. In particular, for an MP enviroment, we must treat $J$ as a sparse matrix in order to overcome the limitation of main memory. In this case, the dense matrix of $J$ need about 36 MB for 50 decimal digits, and hence $3\times 36 \times 10 = 1.08$\ GB for a 10-stage IRK method. On the other hand, the band matrix of $J$ need be about 0.18 MB, and hence,  $3\times 0.18 \times 10 = 5.4$ MB. This problem is sutable for DP-MP Krylov subspace methods in inner iteration. In this case, we use the left preconditioned DP-MP BiCGSTAB and normal DP-MP BiCGSTAB methods.

\tablename\ \ref{table:bruss1d} shows the result of numerical experiments for (\ref{eqn:bruss1d}). The left preconditioned DP BiCGSTAB methods at (\ref{eqn:simple_iterative_ref_l1}) and (\ref{eqn:simple_iterative_ref_l2}) in the DP-MP type iterative refinement method is needed in order to be converged successfully.
\begin{table}
\begin{center}
\caption{1-dimensional Brusselator Problem}\label{table:bruss1d}
\begin{tabular}{|c|c|c|c|c|}\hline
 & \multicolumn{2}{|c|}{Left preconditioned DP BiCGSTAB} & \multicolumn{2}{|c|}{Unpreconditioned} \\ \hline
 & 10 stages & 20 stages & 10 stages & 20 stages \\ \hline
\# steps      & 2966    & 341 & 15056 & 7587  \\ \hline
Comp.Time (s) & 11205   & 3642 & 77717 & 64573 \\ \hline
Max.Rel.Error & $2.2\times 10^{-25}$ & $1.3\times 10^{-20}$& $3.9\times 10^{-25}$ & $3.8\times 10^{-25}$ \\
Min.Rel.Error & $1.3\times 10^{-29}$ & $6.1\times 10^{-23}$& $2.1\times 10^{-27}$ & $2.1\times 10^{-27}$ \\ \hline
\end{tabular}
\end{center}
\end{table}

By comparing the numerical values provided in SPARK3, we confirmed that all elements of the approximations have over true 14 decimal digits. Moreover, we can find that they have about 20 - 29 true decimal digits by comparing the results obtained by using 100 decimal digits MP computation. In the above-mentioned problems, faster computations with the same accuracy of approximations can be achieved with higher-order formulas.

\begin{figure}[htpb]
\begin{center}
\includegraphics[width=.45\textwidth]{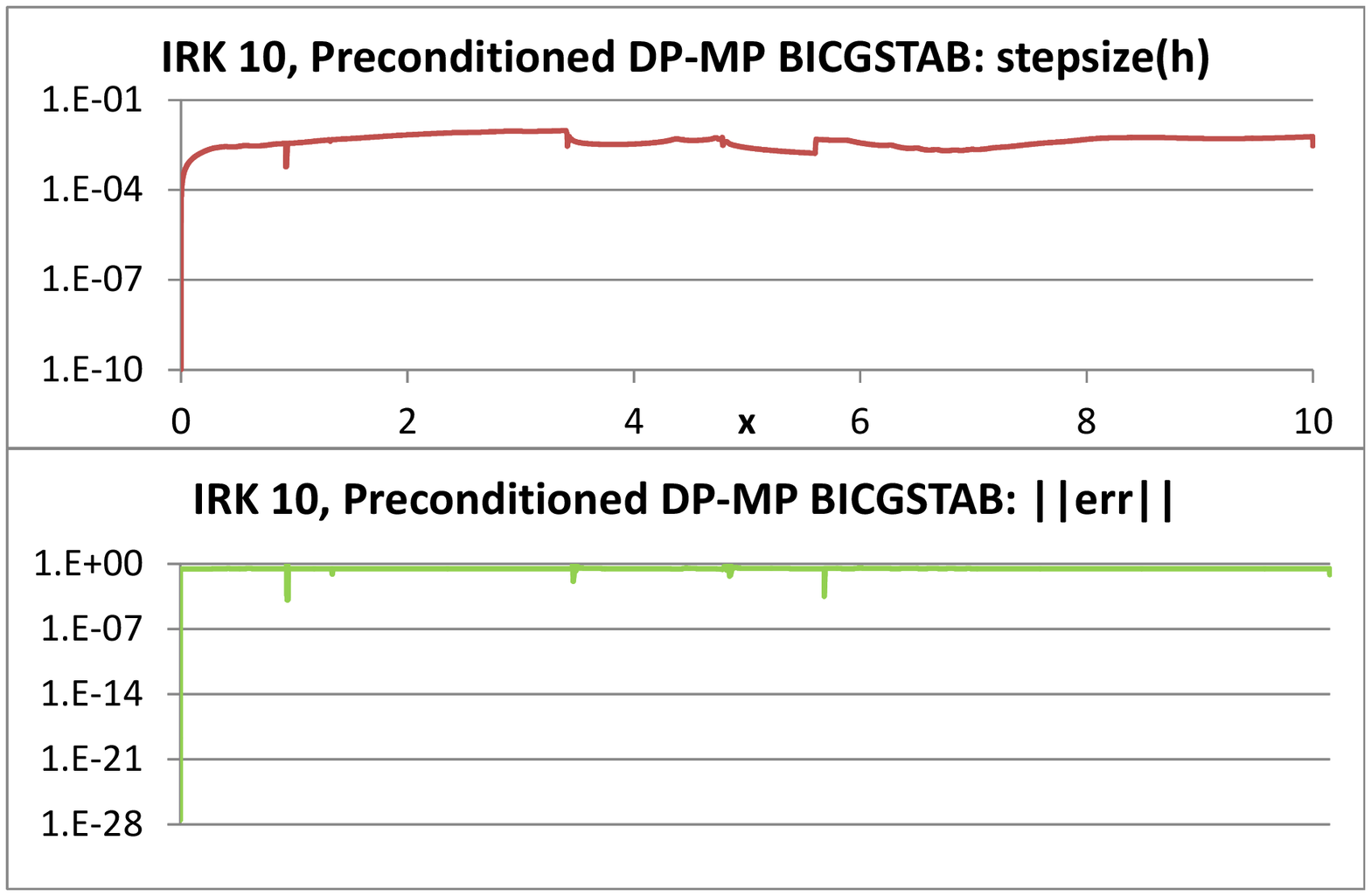}
\includegraphics[width=.45\textwidth]{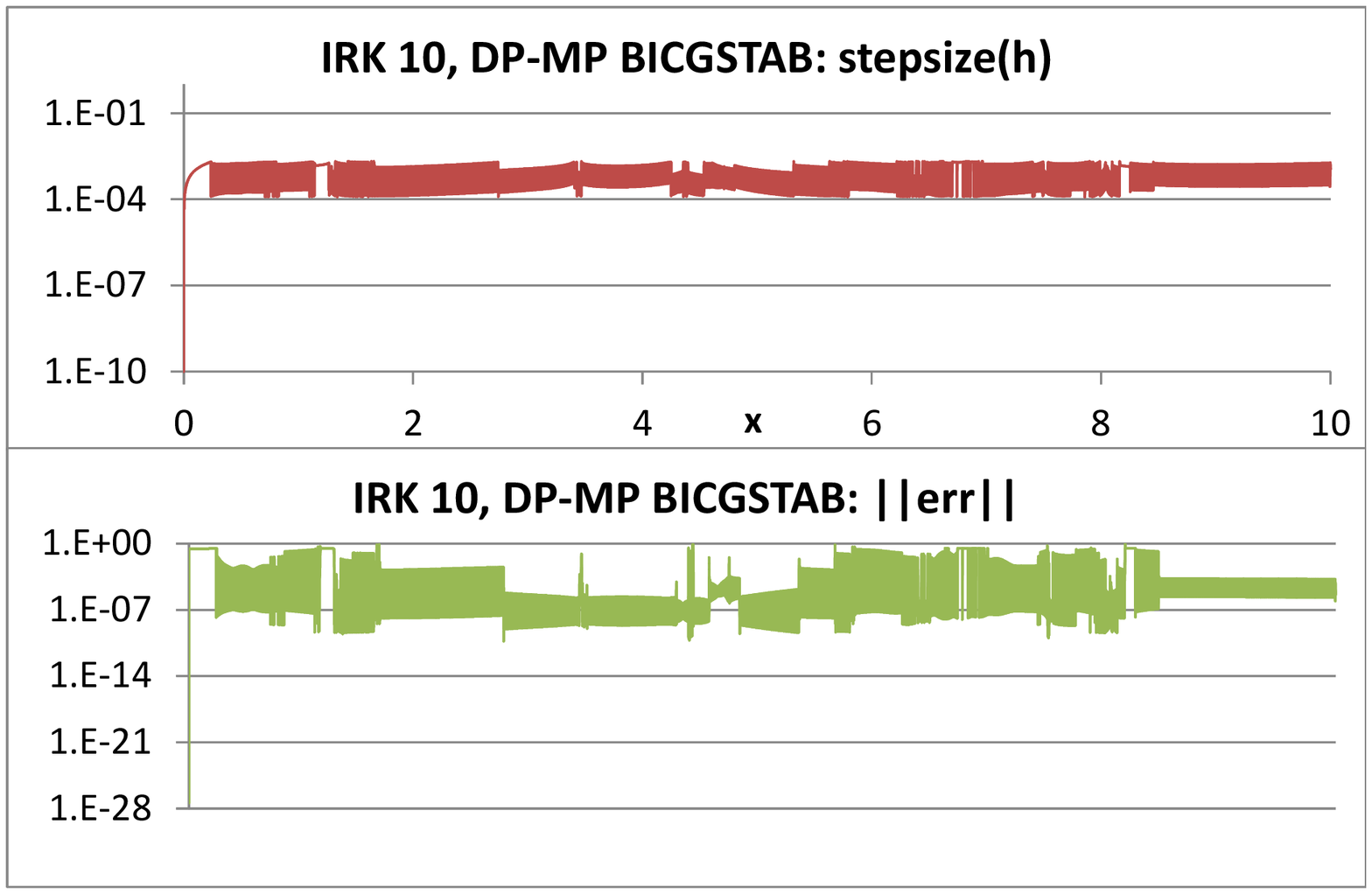}
\caption{History of step size and $\|\mathbf{err}_k\|$ : 10 stages, 20th order : Left preconditioned DP-MP BiCGSTAB(Left), DP-MP BiCGSTAB(Right)}\label{fig:bruss1d_local_error_est}
\end{center}
\end{figure}

However, there are some problems to be solved in order to speed up our MP ODE solver as shown \figurename\ \ref{fig:bruss1d_local_error_est}.

The two sets of graphs of $\|\mathbf{err}_k\|$ and step size show that the right ones are widely wiggled because of non-convergence of the unpreconditioned DP BiCGSTAB method used in the iterative refinement method. If the DP BiCGSTAB method is not convergent, the approximation computation is rejected and a new one is recomputed after the step size becomes smaller at the discretized point. These bottlenecks can be overcome by using predconditionings or other robust DP iterative methods for a system of linear equations in inner iteration.

%
\section{Conlusion and Future Works}

We showed that our MP ODE solver based on high-order IRK formulas can obtain accurate approximations in some problems. However, some problems remain such as the unsupported DP robust linear solver for sparse Jacobi matrix. 

Our final objective is to provide practical and high-performance DP and MP ODE solvers based on high-order IRK methods. 

In order to achieve this objective, we plan to tackle the following issues in the future:
\begin{enumerate}
\item Parallelization of inner iteration in IRK method for multi-core CPU and GPU. For this purpose, we plan to use a well-tuned linear computation library based on LAPACK and BLAS.
\item We plan to accumulate many numerical experiments, especially for ill-conditioned problems requiring MP arithmetic, and we also plan to provide some selection for solvers of system of linear equations in inner iteration in order to optimize the computational time and user-required accuracy.
\end{enumerate}

\section*{Acknowledgments}
I thank Hideko Nagasaka and Masatsugu Tanaka for encouraging me when writing my doctoral thesis that is one of origins of this paper.



\end{document}